\documentclass{article}
\usepackage{amsmath}
\usepackage{amsfonts}
\usepackage{amsthm}
\newcommand{\C}{\mathbb{C}}
\newtheorem{theorem}{Theorem}
\newtheorem{lemma}[theorem]{Lemma}
\newtheorem{cor}[theorem]{Corollary}
\theoremstyle{definition}
\newtheorem{defn}[theorem]{Definition}
\newtheorem{example}[theorem]{Example}
\theoremstyle{remark}
\newtheorem*{acn}{Acknowledgments}
\numberwithin{equation}{section} \numberwithin{theorem}{section}
\DeclareMathOperator{\id}{id} \DeclareMathOperator{\Ad}{Ad}
\DeclareMathOperator{\spa}{span} \DeclareMathOperator{\Iso}{Iso}
\DeclareMathOperator{\cha}{char} \DeclareMathOperator{\Stab}{Stab}
\begin{document}
\title{The Orbit Method for Finite  Groups \\ of Nilpotency Class Two of Odd Order}
\author{Aleksandrs Mihailovs\\
Department of Mathematics\\ Shepherd College\\ Shepherdstown, WV
25443\\ Alec@Mihailovs.com\\ http://www.mihailovs.com/Alec/ }
\maketitle
\begin{abstract}
First, I construct an isomorphism between the categories of
(topological) groups of nilpotency class 2 with 2-divisible center
and  (topological) Lie rings of nilpotency class 2 with
2-divisible center. That isomorphism allows us to construct
adjoint and coadjoint representations as usual. For a finite group
$G$ of nilpotency class 2 of odd order, I construct a basis in its
group algebra $\C[G]$, parameterized by elements of
$\mathfrak{g}^*$  so that the elements of coadjoint orbits form
bases of simple two-side ideals of $\C[G]$. That construction
gives us a one-to-one correspondence between $G$-orbits in
$\mathfrak{g}^*$ and classes of equivalence of irreducible unitary
representations of $G$, implying a very simple character formula.
The properties of that correspondence are similar to the
properties of the analogous correspondence given by Kirillov's
orbit method for nilpotent connected and simply connected Lie
groups. The diagram method introduced in my article \cite{2} and
my thesis \cite{1}, gives us a convenient way to study normal
forms on the orbits and corresponding representations.
\end{abstract}
\setlength{\baselineskip}{1.5\baselineskip}

\section{Introduction}
Fortunately, Kirillov published recently a survey of merits and
demerits of the orbit method \cite{3}, so I can refer to it and
skip a historical introduction here.

\section{Groups of nilpotency class 2 }

By definition, a group $B$ of nilpotency class 2 is a central
extension of an abelian group, i.e.\ there is an exact sequence
\begin{equation}\label{1}
0 \rightarrow A \rightarrow B \rightarrow C \rightarrow 0
\end{equation}
where $A$ and $C$ are abelian groups and $A$ is the center of $B$.
In other words, elements of the group $B$ of nilpotency class 2
can be written as pairs $b=(a,c)$ with $a\in A$ and $c\in C$ so
that
\begin{equation}\label{2}
(a_1,c_1)\cdot (a_2,c_2)=(a_1+a_2+\psi(c_1,c_2),c_1+c_2).
\end{equation}
The associativity of group operation \eqref{2} is equivalent to
the following identity:
\begin{equation}\label{3}
\psi(c_1,c_2)+\psi(c_1+c_2,c_3)=\psi(c_1,c_2+c_3)+\psi(c_2,c_3)
\end{equation}
which means that $\psi$ is a 2-cocycle, $\psi\in \text{C}^2(C,A)$
supposing that the action of $C$ on $A$ is trivial. Note that
substituting either $c_1=c_2=0$, or $c_2=c_3=0$ in \eqref{3}, we
obtain the identities
\begin{equation}\label{3.1}
\psi(0,c)=\psi(c,0)=\psi(0,0).
\end{equation}
For every such cocycle $\psi$, the element $(-\psi(0,0),0)$ is an
identity of $B$ and
\begin{equation}\label{4}
(a,c)^{-1}=(-a-\psi(c,-c)-\psi(0,0),-c)
\end{equation}
is a left inverse to $(a,c)$, which means that formula \eqref{2}
defines a group structure in set $B=A\times C$ for every cocycle
$\psi\in\text{C}^2(C,A)$. In particular, the left inverse
\eqref{4} is a right inverse of $(a,c)$ as well, thus formula
\eqref{3} implies the identity
\begin{equation}\label{5}
\psi(c,-c)=\psi(-c,c).
\end{equation}
In this construction, both mappings
\begin{align}\label{6}
A\rightarrow A\times\{0\},\quad& a\mapsto (a-\psi(0,0),0);\\
C\rightarrow B/A,\quad& c\mapsto (0,c)A \label{7}
\end{align}
are isomorphisms. Choosing other representatives of cosets than
$(0,c)$ in \eqref{7}, for any function $q:C\rightarrow A$,
\begin{equation}\label{8}
C\rightarrow B/A,\quad c\mapsto (0,c)(q(c)-\psi(0,0),0)A=(q(c),c)A
\end{equation}
we obtain the same group structure in $B$. Renaming
\begin{equation}\label{9}
(a,c)_{\text{new}}=(a+q(c),c),
\end{equation}
we obtain from \eqref{2} that the same group structure in $B$ can
be defined by 2-cocycle
\begin{equation}\label{10}
\psi_{\text{new}}(c_1,c_2)=\psi(c_1,c_2)+q(c_1)+q(c_2)-q(c_1+c_2).
\end{equation}
The difference between these old and new cocycles is equal to the
coboundary of 1-chain $q$, thus the group structure in $B$ is
uniquely determined by 2-cocycles modulo coboundaries of 1-chains,
i.e.\, by elements of $\text{H}^2(C,A)$ with trivial action of $C$
on $A$.

\begin{defn}\label{d0}
We'll call a cocycle $\psi\in \text{C}^2(C,A)$ {\em centered} iff
$\psi(0,0)=0$.
\end{defn}

From \eqref{3.1}, for a centered cocycle $\psi$, for every $c\in
C$,
\begin{equation}\label{11}
\psi(0,c)=\psi(c,0)=\psi(0,0)=0.
\end{equation}
Choosing $q(0)=-\psi(0,0)$ in \eqref{10}, we obtain a centered
cocycle. Comparing \eqref{11} with the formula for the identity in
$B$, we see that we can always choose a centered cocycle $\psi$ in
the same cohomology class so the identity element of $B$ was
$(0,0)$.

\begin{defn}\label{d1} An abelian group $A$ is {\em 2-divisible} if
homomorphism
\begin{equation}\label{12}
A\rightarrow A, \quad a\mapsto a+a
\end{equation}
is an automorphism. For an element $a$ of a 2-divisible abelian
group, we denote $a/2$ the image of $a$ under the automorphism
inverse to \eqref{12}.
\end{defn}

\begin{lemma}\label{l1} A finite abelian group is
2-divisible iff it has an odd order. An abelian $p$-group is
2-divisible iff $p$ is odd.
\end{lemma}
\begin{proof} The kernel of \eqref{12} is the subgroup of elements
of order 2. So abelian 2-groups and abelian finite groups of even
order can't be 2-divisible. Finite abelian groups of odd order
don't have elements of order 2, so homomorphism \eqref{12} is
injective, thus its image contains as many elements as $A$ does,
so it is $A$. For abelian $p$-groups with odd $p$, homomorphism
\eqref{12} is injective, and its restriction on every cyclic
subgroup is an automorphism of that subgroup, so \eqref{12} is
surjective as well.
\end{proof}

\begin{defn}\label{d1.1}
We'll call a cocycle $\psi\in \text{C}^2(C,A)$ {\em equalized} iff
$\psi(c,-c)=0$ for all $c\in C$.
\end{defn}

For an equalized cocycle $\psi$, from \eqref{4}, we have
\begin{equation}\label{12.1}
(a,c)^{-1}=(-a,-c)
\end{equation}
for every element $(a,c)\in B$.

\begin{lemma}\label{l2} For a 2-divisible abelian group $A$
and an abelian group $C$, every element of $\text{H}^2(C,A)$ can
be represented by an equalized cocycle.
\end{lemma}
\begin{proof}
If $A$ is a 2-divisible abelian group, for a centered cocycle
$\psi$, choosing $q(c)=-\psi(c,-c)/2$ in \eqref{10}, and combining
it with \eqref{5}, we obtain an equalized cocycle.
\end{proof}

For any $A$, if $C$ doesn't have elements of order 2, an equalized
cocycle can be constructed from a centered cocycle $\psi$ by
choosing $\{q(c),q(-c)\}=\{-\psi(c,-c),0\}$ with an arbitrary
choice of which of them must be equal 0. Here is the example
showing that it can't be obtained for all cases.

\begin{example}\label{ex1}
Let either $B\simeq\text{D}_8$, a dihedral group of order 8, or
$B\simeq\text{Q}_8$, a quaternionic group of order 8 . In both
cases $A\simeq\text{C}_2$ and $C\simeq\text{C}_2\oplus \text{C}_2$
where $\text{C}_n$ denotes a cyclic group of order $n$. We have
$a=-a$ and $c=-c$ for all elements $a\in A$, $c\in C$. If formula
\eqref{12.1} was true, then all nontrivial elements of $B$ would
have order 2, which is false.
\end{example}

If we have a centered or an equalized cocycle and we want to
change it adding the coboundary of a 1-chain $q$ so that the new
cocycle was still centered or equalized, we need the following
conditions on $q$: $q(0)=0$ for centered cocycles, or
$q(-c)=-q(c)$ for all $c\in C$ for equalized cocycles. Note also
that for an equalized cocycle $\psi$ for any $c_1, c_2\in C$,
\begin{equation}\label{14.2}
-\psi(c_1,c_2)=\psi(-c_2,-c_1),
\end{equation}
that follows from the identity $(b_1\cdot b_2)^{-1}=b_2^{-1}\cdot
b_1^{-1}$.

All what was told from \eqref{1} until now, was applicable to all
central group extensions. We didn't explore the fact that $A$ is
exactly the center of $B$. What do we need for that? We need that
cocycle $\psi$ was non-degenerate.

\begin{defn}\label{d2}
A cocycle $\psi$ is {\em non-degenerate} iff  for every $c_1\in
C$, $c_1\neq 0$ there exist such $c_2\in C$ that
\begin{equation}\label{15}
\psi(c_1,c_2)\neq\psi(c_2,c_1).
\end{equation}
\end{defn}

For a non-degenerate cocycle $\psi$, element $b=(a_1,c_1)\in B$
with $c_1\neq 0$ can't be central since it is not commute with
$(0,c_2)$ with $c_2$ satisfying \eqref{15}. Every coboundary of a
1-chain $q$ is symmetric, so adding it to $\psi$ doesn't change
inequality \eqref{15}. That means that all cocycles representing a
cohomology class from $H^2(C,A)$, are either non-degenerate, in
which case we'll call that class non-degenerate; or degenerate,
then we'll call that class degenerate as well.

\section{Lie rings of nilpotency class 2 }

By a Lie ring I mean an abelian group with a bilinear commutator
$[.,.]$ such that $[a,a]=0$ for all $a$, it implies
$[a,b]=-[b,a]$, and satisfying Jacobi identity. A commutative Lie
ring means zero commutator.

By definition, a Lie ring $\mathfrak{b}$ of nilpotency class 2 is
a central extension of a commutative Lie ring, i.e.\ there is an
exact sequence
\begin{equation}\label{2.1}
0 \rightarrow \mathfrak{a} \rightarrow \mathfrak{b} \rightarrow
\mathfrak{c} \rightarrow 0
\end{equation}
where $\mathfrak{a}$ and $\mathfrak{c}$ are commutative Lie rings
and $\mathfrak{a}$ is the center of $\mathfrak{b}$. In other
words, elements of the Lie ring $\mathfrak{b}$ of nilpotency class
2 can be written as pairs $b=(a,c)$ with $a\in \mathfrak{a}$ and
$c\in \mathfrak{c}$ so that
\begin{equation}\label{2.2}
(a_1,c_1)+(a_2,c_2)=(a_1+a_2+\phi(c_1,c_2),c_1+c_2).
\end{equation}
and
\begin{equation}\label{2.2.1}
[(a_1,c_1), (a_2,c_2)]=[\eta(c_1,c_2)-\phi(0,0),0].
\end{equation}
The associativity and commutativity of the operation in
$\mathfrak{b}$ defined in \eqref{2.2} are equivalent to $\phi$
being a symmetric 2-cocycle $\phi\in\text{C}^2(C,A)$ where $C$ and
$A$ are underlying abelian groups of $\mathfrak{c}$ and
$\mathfrak{a}$, correspondingly, supposing that the action of $C$
on $A$ is trivial. The commutator properties are equivalent to the
fact that $\eta:C\times C\rightarrow A$ is a skew-symmetric
bihomomorphism where skew-symmetric means $\eta(c,c)=0$ for all
$c\in C$ which implies $\eta(c_1,c_2)=-\eta(c_2,c_1)$. In other
words, $\eta$ is a 2-cocycle $\eta\in\text{C}^2(\mathfrak{c}, A)$
with trivial action of $\mathfrak{c}$ on A.

In this construction, both mappings
\begin{align}\label{2.6}
\mathfrak{a}\rightarrow \mathfrak{a}\times\{0\},\quad& a\mapsto
(a-\phi(0,0),0);\\ \mathfrak{c}\rightarrow
\mathfrak{b}/\mathfrak{a},\quad& c\mapsto (0,c)+\mathfrak{a}
\label{2.7}
\end{align}
are Lie ring isomorphisms. Choosing other representatives of
cosets than $(0,c)$ in \eqref{2.7}, for any function
$q:C\rightarrow A$,
\begin{equation}\label{2.8}
\mathfrak{c}\rightarrow \mathfrak{b}/\mathfrak{a},\quad c\mapsto
(0,c)+(q(c)-\phi(0,0),0)+\mathfrak{a}=(q(c),c)+\mathfrak{a}
\end{equation}
we obtain the same Lie ring structure in $\mathfrak{b}$. Renaming
\begin{equation}\label{2.9}
(a,c)_{\text{new}}=(a+q(c),c),
\end{equation}
we obtain from \eqref{2.2} that the same Lie ring structure in
$\mathfrak{b}$ can be defined by 2-cocycles
\begin{equation}\label{2.10}
\phi_{\text{new}}(c_1,c_2)=\phi(c_1,c_2)+q(c_1)+q(c_2)-q(c_1+c_2),
\end{equation}
obtaining from $\phi$ by adding a 1-chain, and $\eta$. Note that
the coboundary of a 1-chain $r\in \text{C}^1(\mathfrak{c},A)$ with
the trivial action of $\mathfrak{c}$ on $A$, is
$r(c_1,c_2)=r([c_1,c_2])=0$. Thus the Lie ring structure in
$\mathfrak{b}$ is uniquely determined by elements of
$\text{H}^2_{\text{sym}}(C,A)\oplus \text{H}^2(\mathfrak{c},A)$
with trivial actions of $C$ and $\mathfrak{c}$ on $A$, where
$\text{H}^2_{\text{sym}}(C,A)$ denotes the subgroup of
2-cohomology classes defined by symmetric cocycles. The same as in
the previous section, choosing $q(c)=-\phi(0,0)$ in \eqref{2.10},
we get $\phi_{\text{new}}(0,0)=0$, and $(0,0)=0\in\mathfrak{b}$.

All what was told from \eqref{2.1} until now, was applicable to
all central Lie ring extensions. We didn't explore the fact that
$\mathfrak{a}$ is exactly the center of $\mathfrak{b}$. What do we
need for that? We need that cocycle $\eta$ was non-degenerate.

\begin{defn}\label{d2.2}
A skew-symmetric bihomomorphism $\eta:C\times C\rightarrow A$ is
{\em non-degenerate} iff for every $c_1\in C$, $c_1\neq 0$ there
exist such $c_2\in C$ that
\begin{equation}\label{2.15}
\eta(c_1,c_2)\neq 0.
\end{equation}
\end{defn}

For a non-degenerate cocycle $\eta$, element $b=(a_1,c_1)\in
\mathfrak{b}$ with $c_1\neq 0$ can't be central since its
commutator with $(0,c_2)$ with $c_2$ satisfying \eqref{2.15}, is
not 0.

\section{Lie correspondence}

\begin{theorem}\label{t1}
For a 2-divisible (topological) abelian group $A$ and an abelian
(topological) group $C$ acting trivially on $A$, the following
mapping:
\begin{equation}\label{Lie.1}{\normalfont
L_{\text{c}}:\text{C}^2(C,A)\rightarrow
\text{C}^2_{\text{sym}}(C,A)\oplus \text{C}^2(\mathfrak{c},A),
\quad \psi \mapsto (\phi, \eta)}
\end{equation}
where $\mathfrak{c}$ is the commutative (topological) Lie ring,
underlying abelian group of which is $C$, acting trivially on $A$,
{\normalfont $\text{C}_{\text{sym}}$} denotes symmetric cocycles,
and
\begin{gather}\label{Lie.2}
\phi(c_1,c_2)=\frac{\psi(c_1,c_2)+\psi(c_2,c_1)}{2}, \\
\eta(c_1,c_2)=\psi(c_1,c_2)-\psi(c_2,c_1)\label{Lie.3},
\end{gather}
is an isomorphism, factor of which by coboundaries of 1-chains is
an isomorphism of the cohomology groups
\begin{equation}\label{Lie.4}{\normalfont
L_{\text{h}}: \text{H}^2(C,A)\rightarrow
\text{H}^2_{\text{sym}}(C,A)\oplus \text{H}^2(\mathfrak{c},A).}
\end{equation}
For a 2-divisible (topological) commutative Lie ring
$\mathfrak{a}$ underlying abelian group of which is $A$, and a
commutative (topological) Lie ring $\mathfrak{c}$ acting trivially
on $A$, the following mapping:
\begin{equation}\label{Lie.12}{\normalfont
E_{\text{c}}:\text{C}^2_{\text{sym}}(C,A)\oplus
\text{C}^2(\mathfrak{c},A)\rightarrow \text{C}^2(C,A), \quad
(\phi, \eta)\mapsto \psi}
\end{equation}
where
\begin{equation}\label{Lie.5}
\psi(c_1,c_2)=\phi(c_1,c_2)+\frac{\eta(c_1,c_2)}{2},
\end{equation}
is an isomorphism, factor of which by coboundaries of 1-chains is
an isomorphism of the cohomology groups
\begin{equation}\label{Lie.14}{\normalfont
E_{\text{h}}:\text{H}^2_{\text{sym}}(C,A)\oplus
\text{H}^2(\mathfrak{c},A)\rightarrow \text{H}^2(C,A) .}
\end{equation}
The isomorphisms $L_{\text{c}}$ and $E_{\text{c}}$ are mutually
inverse and the isomorphisms $L_{\text{h}}$ and $E_{\text{h}}$ are
mutually inverse. Cocycle $\phi$ is centered or equalized iff
$\psi$ is centered or equalized, correspondingly. Cocycle $\eta$
is non-degenerate iff $\psi$ is non-degenerate.
\end{theorem}
\begin{proof}
From \eqref{Lie.2}, $\phi$ is a cocycle, by linearity, and it is
symmetric,
\begin{equation}\label{Lie.6}
\phi(c_1,c_2)=\phi(c_2,c_1).
\end{equation}
$\eta$ is skew-symmetric by \eqref{Lie.3}. To check that it is a
bihomomorphism, add the following cocycle identities:
\begin{align}
\psi(c_1+c_2,c_3)+\psi(c_1,c_2)&=\psi(c_1,c_2+c_3)+\psi(c_2,c_3) ,
\label{Lie.7}\\
\psi(c_3+c_1,c_2)+\psi(c_3,c_1)&=\psi(c_3,c_1+c_2)+\psi(c_1,c_2) ,
\label{Lie.8}\\
-\psi(c_1+c_3,c_2)-\psi(c_1,c_3)&=-\psi(c_1,c_3+c_2)-\psi(c_3,c_2).
\label{Lie.9}
\end{align}
After cancelling equal items, we get
\begin{equation}\label{Lie.10}
\psi(c_1+c_2,c_3)+\psi(c_3,c_1)-\psi(c_1,c_3)=\psi(c_3,c_1+c_2)+\psi(c_2,c_3)-\psi(c_3,c_2),
\end{equation}
or
\begin{equation}\label{Lie.11}
\eta(c_1+c_2,c_3)=\eta(c_1,c_3)+\eta(c_2,c_3).
\end{equation}
$L_{\text{c}}$ is a homomorphism by linearity of \eqref{Lie.2} and
\eqref{Lie.3}. Adding coboundary to $\psi$ adds the same
coboundary to $\phi$, and inverse from \eqref{Lie.5}, and we don't
have to worry about adding coboundaries to $\eta$ since
$\text{H}^2(\mathfrak{c},A)= \text{C}^2(\mathfrak{c},A)$, so
$L_{\text{c}}$ defines $L_{\text{h}}$. From the other side, $\psi$
defined in \eqref{Lie.5} is a cocycle by linearity (note, that
$\eta$ is a cocycle since every bihomomorphism from $C\times C$ to
$A$ is an element of  $\text{C}^2(C,A)$). $E_{\text{c}}$ is a
homomorphism by linearity of \eqref{Lie.5}. It is easy to check
that $L_{\text{c}}\circ E_{\text{c}}=\id$ and $E_{\text{c}}\circ
L_{\text{c}}=\id$. Since $L_{\text{h}}$ and $E_{\text{h}}$ are
defined from $L_{\text{c}}$ and $E_{\text{c}}$ by a factorization
by the same coboundaries, they are mutually inverse isomorphisms
as well. The last statement of the theorem immediately follows
from the definitions.
\end{proof}

For a group $B$ of nilpotency class 2 with 2-divisible center
denote $L(B)$ the Lie ring of nilpotency class 2 defined on the
underlying set of $B$ by cocycles \eqref{Lie.2} and \eqref{Lie.3}.
Theorem \ref{t1} tells that this construction doesn't depend on
the choice of the cocycle $\psi$ defining a group structure in
$B$. For a Lie ring $\mathfrak{b}$ of nilpotency class 2 with
2-divisible center denote $E(\mathfrak{b})$ the group of
nilpotency class 2 defined on the underlying set of $\mathfrak{b}$
by a cocycle \eqref{Lie.5}. Theorem \ref{t1} tells that this
construction doesn't depend on the choice of the cocycles $\phi$
and $\eta$ defining a Lie ring structure in $\mathfrak{b}$.

For a homomorphism $f: B_1\rightarrow B_2$ of groups $B_1$ and
$B_2$ of nilpotency class 2 with 2-divisible centers denote
$L(f)=f$ considered as a function from $L(B_1)$ to $L(B_2)$.
Analogously, for a homomorphism $f: \mathfrak{b}_1\rightarrow
\mathfrak{b}_2$ denote $E(f)=f$ considered as a function from
$E(\mathfrak{b}_1)$ to $E(\mathfrak{b}_1)$.

\begin{theorem}\label{t2}
$L$ and $E$ are mutually inverse functors defining an isomorphism
between categories of groups of nilpotency class 2 with
2-divisible center and Lie rings of nilpotency class 2 with
2-divisible center.
\end{theorem}
\begin{proof} By Theorem \ref{t1}, L and E are mutually inverse if
they are functors. All categorical properties would follow
immediately from definitions if we showed that $L(f)$ is a Lie
ring homomorphism for every group homomorphism $f$ and $E(f)$ is a
group homomorphism for every Lie ring homomorphism $f$.

By Lemma \ref{l2}, we can choose equalized cocycles $\psi_1$ and
$\psi_2$ defining group structures in $B_1$ and $B_2$. Then
symmetric cocycles $\phi_1$ and $\phi_2$ defined by Theorem
\ref{t1}, are equalized as well. Multiplying the left hand sides
and the right hand sides of the following identities:
\begin{multline}\label{Lie.15}
f(a_1+a_2+\psi_1(c_1,c_2),c_1+c_2)=f(a_1,c_1)\cdot f(a_2,c_2)\\
=(\alpha_1,\gamma_1)\cdot
(\alpha_2,\gamma_2)=(\alpha_1+\alpha_2+\psi_2(\gamma_1,\gamma_2)
,\gamma_1+\gamma_2) ,
\end{multline}
\begin{multline}
f(-a_1-a_2-\psi_1(c_2,c_1),-c_1-c_2)=(f(a_2,c_2)\cdot
f(a_1,c_1))^{-1}\\ =((\alpha_2,\gamma_2)\cdot
(\alpha_1,\gamma_1))^{-1}=(-\alpha_1-\alpha_2-\psi_2(\gamma_2,\gamma_1),
-\gamma_1-\gamma_2) \label{Lie.16}
\end{multline}
we get
\begin{equation}\label{Lie.17}
f(\eta_1(c_1,c_2),0)=(\eta_2(\gamma_1,\gamma_2),0) .
\end{equation}
Dividing the central elements by 2 and inverting, we get
\begin{equation}\label{Lie.18}
f(-\frac{\eta_1(c_1,c_2)}{2},0)=
(-\frac{\eta_2(\gamma_1,\gamma_2)}{2},0) .
\end{equation}
Multiplying \eqref{Lie.18} and \eqref{Lie.15}, we get
\begin{equation}\label{Lie.19}
f(a_1+a_2+\phi_1(c_1,c_2),c_1+c_2)=(\alpha_1+\alpha_2+\phi_2(\gamma_1,\gamma_2)
,\gamma_1+\gamma_2).
\end{equation}
So
\begin{equation}\label{Lie.20}
L(f)(b_1+b_2)=L(f)(b_1)+L(f)(b_2)
\end{equation}
for $b_1=(a_1,c_1)$, $b_2=(a_2,c_2)$. This formula together with
\eqref{Lie.17} rewritten as
\begin{equation}\label{Lie.21}
L(f)([b_1,b_2])=[L(f)(b_1), L(f)(b_2)]
\end{equation}
means that $L(f)$ is a Lie ring homomorphism.

Similarly, for a Lie ring homomorphism $f$, multiplying
\eqref{Lie.19} and the identity obtained from \eqref{Lie.18} by
changing $-$ to $+$, we get \eqref{Lie.15} which means that $E(f)$
is a group homomorphism.
\end{proof}

\begin{lemma}\label{Lie.l1}
For any elements $b, b_1, b_2$ of a group $B$ of nilpotency class
2 with 2-divisible center,
\begin{gather}\label{Lie.22}
b^{-1}=-b ,\\ \label{Lie.23} b_1\cdot
b_2=b_1+b_2+\frac{[b_1,b_2]}{2} ,\\ \label{Lie.24} b_1\cdot
b_2\cdot b_1^{-1}=b_2+[b_1,b_2] ,\\ \label{Lie.25} b_1\cdot
b_2\cdot b_1^{-1}\cdot b_2^{-1}=[b_1,b_2] .
\end{gather}
If $b_1$ and $b_2$ commute, then
\begin{equation}\label{Lie.26}
b_1\cdot b_2=b_1+b_2 .
\end{equation}
\end{lemma}
\begin{proof}
Formula \eqref{Lie.22} is true because by Lemma \ref{l2} we can
choose an equalized cocycle $\psi$ defining group structure in $B$
and it corresponds by Theorem \ref{t1} to the equalized cocycle
$\phi$ defining additive group structure in $L(B)$. Formula
\eqref{Lie.26} is true because $\psi(c_1,c_2)=\phi(c_1,c_2)$ for
commuting $b_1=(a_1,c_1)$ and $b_2=(a_2,c_2)$. We already used
\eqref{Lie.23} at the end of the proof of Theorem \ref{t2}.
Formula \eqref{Lie.25} telling that the group commutator in $B$
coincides with the Lie ring commutator in $L(B)$ is true since the
product of left hand sides of formulas \eqref{Lie.15} and
\eqref{Lie.16} equals left hand side of \eqref{Lie.17}. Formula
\eqref{Lie.24} can be obtained by right multiplication of both
sides of \eqref{Lie.25} by $b_2$ and using \eqref{Lie.26}.
\end{proof}

\begin{defn}\label{Lie.d1}
A group $B$ is {\em 2-rootable} if mapping
\begin{equation}\label{12.3}
B\rightarrow B, \quad b\mapsto b^2
\end{equation}
is a bijection.
\end{defn}

\begin{cor}\label{Lie.l2}
A finite group of nilpotency class 2 is 2-rootable iff it has an
odd order. A $p$-group of nilpotency class 2 is 2-rootable iff $p$
is odd.
\end{cor}
\begin{proof}
Groups of even order and 2-groups have elements of order 2, so the
identity covers more than once by \eqref{12.3}, and these groups
can't be 2-rootable. Finite groups of odd order and $p$-groups
with odd $p$ have 2-divisible center, by Lemma \ref{l1}, so
$b^2=b+b$ by \eqref{Lie.26}, and Lemma \ref{l1} applied to $L(B)$,
completes the proof.
\end{proof}

\section{The orbit method}

\begin{lemma}\label{Orb.l1}
For a group $B$ of nilpotency class 2 with 2-divisible center,
formula
\begin{equation}\label{Orb.1}
\Ad(b)(l)=L(b\cdot l\cdot b^{-1})
\end{equation}
defines a structure of left $B$-module in the underlying abelian
group of $L(B)$.
\end{lemma}
\begin{proof}
Conjugation by $b\in B$ is an automorphism of $B$, so $\Ad(b)$ is
a Lie ring automorphism by Theorem \ref{t2}. The formula
\begin{equation}\label{Orb.2}
\Ad(b_1)\circ \Ad(b_2)=\Ad(b_1\cdot b_2)
\end{equation}
follows directly from the definition \eqref{Orb.1}.
\end{proof}

As usual, we'll call $\Ad$ adjoint representation of $B$. Denote
$L(B)^*$ the group of (unitary) characters of the underlying
abelian group of $L(B)$. Define coadjoint representation $\Ad^*$
as the dual to $\Ad$ left action of $B$ in $L(B)^*$:
\begin{equation}\label{Orb.3}
\Ad^*(b)(\chi)(l)=\chi(\Ad(b^{-1})(l)).
\end{equation}

\begin{lemma}\label{Orb.l2}
For a group $B$ of nilpotency class 2 with 2-divisible center, for
every $b\in B$, $l\in L(B)$, $\chi\in L(B)^*$,
\begin{gather}\label{Orb.4}
\Ad(b)(l)=l+[b,l] ,\\ \Ad^*(b)(\chi)(l)=\chi(l-[b,l])
.\label{Orb.5}
\end{gather}
\end{lemma}
\begin{proof}
Both formulas follow from \eqref{Lie.24} and definitions
\eqref{Orb.1} and \eqref{Orb.3}.
\end{proof}

Denote $\mathcal{O}(B)$ the set of orbits of the coadjoint
representation of $B$.

\begin{theorem}[Orbit method]\label{Orb.t1}
For a finite group $B$ of nilpotency class 2 of odd order,
elements
\begin{equation}\label{Orb.6}
\left( X_{\chi}=\sum_{l\in L(B)} \chi(l)\; l\right)_{\chi\in
L(B)^*}
\end{equation}
form an orthonormal basis in $\C[B]$. For every orbit $\Omega$ of
the coadjoint representation of $B$, the subspace
\begin{equation}\label{Orb.7}
V_{\Omega} = \spa(X_{\chi})_{\chi\in \Omega}
\end{equation}
is a two-side ideal of a group algebra $\C[B]$, the restriction of
the regular representation of $B$ in $\C[B]$ on $V_{\Omega}$ is
isotypic, and
\begin{equation}\label{Orb.8}
\C[B]=\bigoplus_{\Omega\in \mathcal{O}(B)} V_{\Omega}
\end{equation}
is the decomposition of the regular representation of $B$ in
$\C[B]$ in the direct sum of isotypic components.
\end{theorem}
\begin{proof}
Elements \eqref{Orb.6} form an orthonormal basis in $\C[L(B)]$ and
dot products in $\C[L(B)]$ and $\C[B]$ are the same.
\begin{multline}\label{Orb.9}
b\; X_{\chi}=\sum_{l\in L(B)} \chi(l)\; b l=\sum_{g\in B}
\chi(b^{-1} g) g =\sum_{g\in
B}\chi\left(-b+g+\frac{[-b,g]}{2}\right) g\\ = \chi(-b)\sum_{g\in
B}\Ad^*\left(\frac{b}{2}\right)(\chi) (g)\; g=\chi(-b)\;
X_{\Ad^*(b/2)(\chi)}\in V_{\Omega}.
\end{multline}
Similarly,
\begin{multline}\label{Orb.10}
X_{\chi} b=\sum_{l\in L(B)} \chi(l)\; l b =\sum_{g\in B} \chi( g
b^{-1}) g =\sum_{g\in B}\chi\left(g-b+\frac{[g,-b]}{2}\right) g\\
= \chi(-b)\sum_{g\in B}\Ad^*\left(\frac{-b}{2}\right)(\chi) (g)\;
g=\chi(-b)\; X_{\Ad^*(-b/2)(\chi)}\in V_{\Omega}.
\end{multline}
So $V_{\Omega}$ is a two-side ideal of $\C[B]$. Because $X_{\chi}$
are orthogonal for different $\chi$, these ideals $V_{\Omega}$ are
orthogonal for different $\Omega$. Since $\C[B]$ is a semi-simple
associative algebra and every component in the direct sum of the
isotypic components
\begin{equation}\label{Orb.11}
\C[B]=\bigoplus_{\tau\in \Hat{B}} \Iso(\tau)
\end{equation}
is simple, every $V_{\Omega}$ is either one of these components,
or a direct sum of a few of them. But the number of coadjoint
orbits coincides with the number of conjugate classes of $B$, by
Duality Lemma \ref{Orb.l3}, i.e.\ with the number of irreducible
unitary representations of $B$, i.e.\ with the number of
components in the sum \eqref{Orb.11}. Thus every $V_{\Omega}$
coincides with isotypic component $\Iso (\tau)\simeq n \tau$ for
an irreducible unitary representation $\tau$ where $n=\dim\tau$,
and different orbits correspond to different irreducible
representations.
\end{proof}

\begin{lemma}[Duality Lemma]\label{Orb.l3}
Let a finite abelian group $M$ be a left $G$-module for a finite
group $G$. Denote $M^*$ the group of unitary characters of $M$ and
define the structure of a dual left $G$-module in $M^*$ by formula
\begin{equation}\label{Orb.12}
g\chi (m)=\chi (g^{-1}m) .
\end{equation}
Then the number of $G$-orbits in $M^*$ is the same as the number
of $G$-orbit in $M$.
\end{lemma}
\begin{proof}
Denote $T_M$ the representation of $G$ in $\C[M]$ defined by
formula $T_M(g)m=gm$. If element
\begin{equation}\label{Orb.14}
x=\sum_{m\in M}x_m m
\end{equation}
is an invariant of $T_M$, then for every orbit $\Xi$ of $G$ in $M$
all coefficients $x_m$ with $m\in \Xi$ are the same. That means
that elements
\begin{equation}\label{Orb.15}
e_{\Xi}=\sum_{m\in \Xi} m
\end{equation}
form a basis in the space of invariants of $T_M$, and the
dimension of the space of invariants equals number of $G$-orbits
in $M$. The same is true for $M^*$: the number of $G$-orbits in
$M^*$ equals the dimension of the space of invariants of
representation $T_{M^*}$ defined by formula
$T_{M^*}(g)\chi=g\chi$. By construction, the space $\C[M^*]$ is
the dual space to $\C[M]$ and representations $T_M$ and $T_{M^*}$
are dual. Thus if
\begin{equation}\label{Orb.16}
T_M\simeq \bigoplus_{\tau\in \Hat{G}}n_{\tau} \tau
\end{equation}
is the decomposition of $T_M$ in the sum of irreducible
representations, then
\begin{equation}\label{Orb.17}
T_{M^*}\simeq \bigoplus_{\tau\in \Hat{G}}n_{\tau} \tau^*
\end{equation}
is the decomposition of $T_{M^*}$ in the sum of irreducible
representations. The dimensions of the spaces of the invariants
equal to multiplicities of the trivial representation in
\eqref{Orb.16} and \eqref{Orb.17}, which are the same since the
trivial representation is self-dual.
\end{proof}

Formulas \eqref{Orb.8} and \eqref{Orb.11} set two different
isomorphisms between the set of coadjoint orbits $\mathcal{O}$ and
the set of classes of equivalency of irreducible unitary
representations $\Hat{B}$, depending on what regular
representation of $B$ in $C(B)$ we use.
\begin{defn}\label{Orb.d1}
Denote $\tau(\Omega)$ and $\Omega(\tau)$ the irreducible unitary
representation class and coadjoint orbit so that
\begin{equation}\label{Orb.19}
V_{\Omega(\tau)}=\Iso(\tau(\Omega))\simeq n \tau(\Omega)
\end{equation}
where $n=\dim \tau (\Omega)$ and the regular representation of $B$
in $\C[B]$ used in \eqref{Orb.11} and \eqref{Orb.19}, is defined
as
\begin{equation}\label{Orb.19.1}
R(g)x=xg^{-1}.
\end{equation}
\end{defn}

We need to use the regular representation given by
\eqref{Orb.19.1} to ensure for abelian $B$ for an orbit containing
one character, correspondence to that character.

\begin{cor}[Dimension formula] \label{Orb.c1}
For a finite group $B$ of nilpotency class 2 of odd order, every
coadjoint orbit $\Omega$ of $B$ has $n^2$ elements where
$n=\dim\tau(\Omega)$. In other words,
\begin{equation}\label{Orb.18}
\dim \tau (\Omega)=\sqrt{\# \Omega (\tau)} .
\end{equation}
\end{cor}
\begin{proof}
From \eqref{Orb.7},
\begin{equation}\label{Orb.20}
\dim V_{\Omega}= \# \Omega.
\end{equation}
Formula \eqref{Orb.18} follows directly from here and
\eqref{Orb.19}
\end{proof}

\begin{lemma}[Stabilizer lemma]\label{Orb.l4}
For a group $B$ of nilpotency class 2 with 2-divisible center, for
any coadjoint orbit $\Omega$ and $\chi_1,\chi_2\in\Omega$
\begin{equation}\label{Orb.21}
\Stab \chi_1=\Stab \chi_2
\end{equation}
and
\begin{equation}\label{Orb.22}
\chi_1(b)=\chi_2(b)
\end{equation}
for any $b\in \Stab \chi_1=\Stab \chi_2$ where $\Stab \chi$
denotes the stabilizer of $\chi$.
\end{lemma}
\begin{proof}
If $b\in \Stab \chi_1$, then for all $l\in L(B)$
\begin{equation}\label{Orb.24}
\chi_1(l)=\chi_1(l-[b,l])=\chi_1(l)/\chi_1([b,l])
\end{equation}
so $\chi_1([b,l])=1$. Since $\chi_1$ and $\chi_2$ are on the same
orbit, there is such $g\in B$ that for all $l\in L(B)$
\begin{equation}\label{Orb.23}
\chi_2(l)=\chi_1(l-[g,l]) .
\end{equation}
Then
\begin{equation}\label{Orb.25}
\chi_2([b,l])=\chi_1([b,l]-[g,[b,l]])=\chi_1([b,l])=1
\end{equation}
so
\begin{equation}\label{Orb.27}
\chi_2(l-[b,l])=\chi_2(l)/\chi_2([b,l])=\chi_2(l)
\end{equation}
so $b\in \Stab \chi_2$ and we proved \eqref{Orb.21}. Now
\begin{equation}\label{Orb.26}
\chi_2(b)=\chi_1(b-[g,b])=\chi_1(b)\chi_1([b,g])=\chi_1(b)
\end{equation}
\end{proof}

\begin{theorem}[Character formula]\label{Orb.t2}
For a finite group $B$ of nilpotency class 2 of odd order, for any
$b\in B$ and $\chi\in\Omega$
\begin{equation}\label{Orb.28}
\cha \tau (\Omega) (b) =
  \begin{cases}
    n\chi(b) & \text{if $b\in \Stab\chi$}, \\
    0 & \text{otherwise},
  \end{cases}
\end{equation}
where $n=\dim \tau (\Omega)$.
\end{theorem}

\begin{proof}
Find the character of the restriction of regular representation
\eqref{Orb.19.1} to $V_{\Omega}$. From \eqref{Orb.10}, if $b/2\in
\Stab \chi$, then from Stabilizer Lemma \ref{Orb.l4}, $b$ acts in
$V_{\Omega}$ by scalar multiplication on $\chi(b)$, so
\begin{equation}\label{Orb.29}
\cha n\tau(\Omega)(b)=\dim V_{\Omega}\chi(b)=n^2\chi(b)
\end{equation}
Note that the cyclic subgroup generated by $b$ is of odd order, so
it contains $b/2$, that means that $b$ and $b/2$ either belong to
$\Stab \chi$, or not, simultaneously. Now, if $b\not\in \Stab
\chi$, then again from \eqref{Orb.10} and from Stabilizer Lemma
\ref{Orb.14}, all diagonal elements of the matrix of the action of
$b$ in $V_{\Omega}$ are zeroes, so
\begin{equation}\label{Orb.30}
\cha n\tau(\Omega)(b)=0
\end{equation}
Dividing \eqref{Orb.29} and \eqref{Orb.30} by $n$, we get
\eqref{Orb.28}.
\end{proof}

\begin{cor}\label{Orb.c4}
For a finite abelian group $B$ of odd order, for
$\Omega=\{\chi\}$,
\begin{equation}\label{Orb.32}
\tau(\Omega)=\chi
\end{equation}
\end{cor}
\begin{proof}
It follows from \eqref{Orb.30} for $n=1$ and $\Stab \chi=B$.
\end{proof}

\begin{cor}\label{Orb.c3}
For a finite group $B$ of nilpotency class 2 of odd order,
\begin{equation}\label{Orb.31}
\tau(-\Omega)=\tau(\Omega)^* .
\end{equation}
\end{cor}
\begin{proof}
Again it follows directly from \eqref{Orb.28}.
\end{proof}

\begin{acn}
I would like to thank Martin Isaacs, Wolfgang Kappe and Peter
Morris for useful discussions.
\end{acn}

\end{document}